\documentclass{amsart}
\usepackage{amsfonts,amssymb,amsmath,amsthm,bbm,xcolor,ulem,graphicx}

\newtheorem{theorem}{Theorem}[section]
\newtheorem{lemma}[theorem]{Lemma}

\newtheorem{proposition}[theorem]{Proposition}
\newtheorem{corollary}[theorem]{Corollary}

\newtheorem{note}[theorem]{Notation}
\newtheorem{remark}[theorem]{Remark}
\newtheorem{conjecture}[theorem]{Conjecture}
\theoremstyle{definition}
\numberwithin{equation}{section}

\def\a{\mathfrak{a}}

\def\R{{\bf R}}

\def\C{{\bf C}}

\def\diag{\mathop{\hbox{diag}}}

\def\to{\rightarrow}

\begin{document}

\allowdisplaybreaks

\begin{center}
{\bf Sharp estimates  for the  hypergeometric functions related to\\  root systems of type $A$ and of rank 1}\\
Piotr Graczyk and Patrice Sawyer
\end{center}

\section*{Abstract}
In this article, we conjecture exact estimates for the Weyl-invariant Opdam-Cherednik hypergeometric functions.  We prove the conjecture for the root system $A_n$
and for all  rank 1 cases.  We provide other evidence that the conjecture might be true in general.

\section{Introduction and Conjecture}

\subsection{Basics on Opdam-Cherednik analysis}
In Opdam-Cherednik analysis, the ``curved'' counterpart of Dunkl analysis for a root system $\Sigma$ on $\R^d$, a crucial role is played by the Opdam-Cherednik kernel $G_k(X,Y)$.
Finding good estimates of  the kernels $G_k$ is therefore important.  In this paper we conjecture exact estimates of the $W$-radial Opdam-Cherednik kernel or hypergeometric functions related to root systems.  We prove these estimates in the case of the root systems of type $A$ and for all rank one cases.

It is interesting to note that our ``guess'' for the behaviour of the hypergeometric functions related to the root system $A_n$ was completely informed by the rank one case $n=1$.  It is therefore encouraging for our conjecture that we were able to verify it for the other root system of rank 1, namely $BC_1$.

For a good introduction on Opdam-Cherednik theory, the reader should consider the paper \cite{Opdam}
and the Lecture Notes \cite{OpdamLN} by Opdam (see also \cite{AAS, Schapira}).  We provide here some details and notations on Opdam-Cherednik analysis.

For every root $\alpha\in \Sigma$, let $\sigma_\alpha(X)=X-2\,\frac{\langle \alpha,X\rangle}{\langle \alpha,\alpha\rangle}\,\alpha$. 
The Weyl group $W$ associated to the root system is generated by the reflection maps $\sigma_\alpha$.

A function $k: \Sigma \to [0,\infty)$ is called a multiplicity function if it is invariant under the action of $W$ on $\Sigma$.

Let $\partial_\xi$ be the derivative in the direction of $\xi\in\R^d$.  The Dunkl-Cherednik or Cherednik operators indexed by $\xi$ are then given by
\begin{align*}
D_\xi\,f(X)&=\partial_\xi\,f(X)+\sum_{\alpha\in \Sigma^+}\,k_\alpha\,\alpha(\xi)\,\frac{f(X)-f(\sigma_\alpha\,X)}{1-e^{-\alpha}}-\rho(k)(\xi)\,f(X),
\end{align*}
where $\rho(k)= \sum_{\alpha\in \Sigma_+} k(\alpha)\alpha$.
The $D_\xi$'s, $\xi\in\R^d$, form a commutative family.

For fixed $ Y\in\R^d$, the  kernel $G_k(\cdot,\cdot)$ is  the only real-analytic solution to the system
\begin{align*}
\left.D_\xi(k)\right\vert_X\,G_k(X,Y)
=\langle \xi,Y\rangle\, G_k(X,Y),~\forall\xi\in\R^d
\end{align*}
with $G_k(0,Y)=1$.  In fact, $G_k$ extends to a holomorphic function  on $(\R^d+i\,U)\times \C^d$ where $U$ is a neighbourhood of 0 (refer to \cite[Th.{} 3.15]{Opdam}).

Its $W$-invariant version $G^W_k(X,\lambda)$ is called  a hypergeometric function.  One notes that in \cite{AAS}, the authors use the term ``hypergeometric function'' for the function  $G_k(X,Y)$.

The (Weyl-invariant) hypergeometric functions related to root systems are the extension of the spherical functions for noncompact symmetric spaces $\phi_\lambda$ to arbitrary positive multiplicities.
In this paper we use the latter terminology and notation.  We have
\begin{align*}
\phi_\lambda(X) =G^W_k(X,\lambda)=\frac{1}{|W|}\sum_{w\in w}\,G_k(w\cdot X,\lambda)
\end{align*}
and $\phi_\lambda(X)$ is the only real-analytic solution of the system
\begin{align*}
\left.p(D_{{\bf e}_1},\dots,D_{{\bf e}_d})(k)\right\vert_X\,
\phi_\lambda(X)
=p(\lambda)\,\phi_\lambda(X), \qquad ~\forall \lambda \in\R^d
\end{align*}
for every Weyl-invariant polynomial $p$ (here ${\bf e}_1$, \dots, ${\bf e}_d$ represent the standard basis on $\R^d$).

Let $\omega_k(X):=
\prod_{\alpha\in \Sigma^+}\,\vert\sinh\langle \alpha,X\rangle\vert^{2\,k(\alpha)}
$
be the Opdam-Cherednik weight function  on $\R^d$.
Recall that the Opdam-Cherednik transform of a $W$-invariant function $f$ on $\R^d$
\begin{align*}
\hat f(\lambda):= c_k^{-1} \int f(x) \phi_{-i\lambda}(X)\omega_k(X)dX,\qquad \lambda\in \R^d,
\end{align*}
plays the role of the spherical Fourier transform in $W$-invariant Opdam-Cherednik analysis.
 
\subsection{Conjecture on sharp bound for the hypergeometric functions related to root systems}

The notation $f\asymp g$ in a domain $D$ means that there exists $C_1>0$ and $C_2>0$ such that
$C_1\,g(x)\leq f(x)\leq C_2\,g(x)$ with $C_1$ and $C_2$ independent of $x\in D$. 

\begin{conjecture}\label{conj:all}
If $X$, $\lambda\in\overline{\a^+}$, then we have for any root system 
\begin{align*}
 \phi_\lambda(e^X)\asymp 
e^{(\lambda-\rho)(X)}\,\prod_{\alpha\in\Sigma^{++}
}\,\frac{1+\alpha(X)}{(1+\alpha(\lambda)\,\alpha(X))^{k(\alpha)+k(2\alpha)} }
\,\left(\frac{1+\alpha(\lambda)\,(1+\alpha(X))}{1+\alpha(\lambda)}\right)^{
k(\alpha)+k(2\alpha)-1}
\end{align*}
where $\Sigma^{++}$ is the set of indivisible positive roots and
$\Sigma^{++}_{\lambda}=\{\alpha\in  \Sigma^{++}\colon \alpha(\lambda)=0\}$.
\end{conjecture}

%%%%%%%%%%%%%%%%%%%%%%%%%%%%%%%%%%%%%%%%%%%%

\begin{remark}
According to Conjecture \ref{conj:all}, we have
\begin{align*}
\phi_\lambda(e^X)\asymp 
e^{(\lambda-\rho)(X)}\,\prod_{\alpha\in\Sigma^{++}
}(1+\alpha(X)) \,f_\alpha(\lambda,X)
\end{align*}
where
the function
\begin{align*}
f_\alpha(\lambda,X)=\frac{1}{(1+\alpha(\lambda)\,\alpha(X))^k }
\,\left(\frac{1+\alpha(\lambda)\,(1+\alpha(X))}{1+\alpha(\lambda)}\right)^{k-1}, \qquad k=k(\alpha)+k(2\,\alpha),
\end{align*}
codifies in one $k$-rational function (i.e.{} a rational function with powers $k>0$  non-necessarily an integer)
four possible power function asymptotics:
\begin{align}
f_\alpha(\lambda,X) &\asymp\left\lbrace
\begin{array}{cl}
1&\hbox{if $\alpha_X\,\alpha_\lambda\leq1$},\\
 1/(\alpha_\lambda\,\alpha_X)^k&\hbox{if $\alpha_X\,\alpha_\lambda\geq1$,$~\alpha_X\leq 1$},\\
1/(\alpha_\lambda\,\alpha_X)&\hbox{if $\alpha_X\,\alpha_\lambda\geq1$,~$\alpha_X\geq 1$,~$\alpha_\lambda\leq1$},\\
1/(\alpha_\lambda^k\,\alpha_X)&\hbox{if $\alpha_X\,\alpha_\lambda\geq1$,~ $\alpha_X\geq 1$,$~\alpha_\lambda\geq1$}.
\end{array}
\right.\label{regions}
\end{align}
\end{remark}

\setlength{\unitlength}{0.8 cm}
\begin{figure}[ht]
\begin{picture}(6,0.0)(0.0,5.0)

\put(0,0){\vector(1,0){5.5}}
\put(0,0){\vector(0,1){5.5}}

\put(4.5,0.4){$\scriptstyle \alpha_\lambda\,\alpha_X=1$}

\qbezier(0.2,5.0)(0.8,1.25)(1.0,1.0)
\qbezier(1.0,1.0)(2.0,0.5)(5.0,0.2)

\put(1.0,1.0){\line(1,0){4}}
\put(1.0,1.0){\line(0,1){4}}

\put(0.5,0.5){\tiny  I}
\put(0.6,3.3){\tiny  II}
\put(3.0,0.7){\tiny III}
\put(2.5,2.5){\tiny IV}

\put(-0.6,5.5){$\scriptstyle \alpha_X$}
\put(5.5,-0.3){$\scriptstyle \alpha_\lambda$}
\end{picture}
{\vskip 4.0cm}
\caption{The four regions\label{R}}
\end{figure}

%%%%%%%%%%%%%%%%%%%%%%%%%%%%%%%%%%%%%%%%%%%%%

The main result of this paper is the proof of Conjecture \ref{conj:all} in the $A_n$ case.

\begin{theorem}\label{trig}
If $X$, $\lambda\in\overline{\a^+}$, then we have for the root systems of  type $A$,
\begin{align*}
 \phi_\lambda(e^X)\asymp 
e^{(\lambda-\rho)(X)}\,\prod_{\alpha\in\Sigma^+
}\,\frac{1+\alpha(X)}{(1+\alpha(\lambda)\,\alpha(X))^k }
\,\left(\frac{1+\alpha(\lambda)\,(1+\alpha(X))}{1+\alpha(\lambda)}\right)^{k-1}. 
\end{align*}
\end{theorem}

Naturally, this theorem is consistent with the results obtained in \cite{PGPS1} for the complex case ($k=1$).

\subsection{Outline of the paper}
In Section \ref{tech}, we introduce notation and some results that will be useful to prove the various upper and lower estimates.
The proof of Theorem \ref{trig} is found in Section \ref{P}.
We conclude with Section \ref{evidence} where we present other evidence for Conjecture \ref{conj:all} including the proof for the root system $BC_1$. In Section \ref{comp}, we show that our conjecture is consistent with some known estimates \cite{Narayanana}. 

\section{Notation and technical results}\label{tech}

\begin{note}\label{notation}
We will write $f(x)\lesssim g(x)$ ($f(x)\gtrsim g(x)$) for $x\in D$ if there exists a constant $C>0$ independent of $x$ such that $f(x)\leq C\,g(x)$ ($f(x)\geq C\,g(x)$) for all $x\in D$.
\end{note}

In what follows, $\Sigma_n^+$ will be the set of positive roots of the root system $A_n$.

We introduce here some technical results.

\begin{lemma}\label{incr}
Assume $a\geq 0$.  Then for $u\geq 0$, the functions $F_1$, $F_2$, $F_3$, $F_4$ and $F_5$ defined by
\begin{align*}
\begin{array}{lll}
F_1(u)=\frac{u}{1+a\,u}& F_2(u)=\frac{u\,(1+a\,(1+u))^{k-1}}{(1+a\,u)^k},~0<k\leq 1, &F_3(u)=\frac{u\,(1+a\,(1+u))}{(1+u)\,(1+a\,u)}\\
F_4(u)=\frac{1+a\,u}{1+a\,(1+u)}&F_5(u)=\frac{1+u}{1+a\,(1+u)}&
\end{array}
\end{align*}
are all increasing functions of $u$.
\end{lemma}

\begin{proof}
It suffices to compute the derivatives $F_1'(u)$, $F_2'(u)$, $F_3'(u)$, $F_4'(u)$, $F_5'(u)$ which are easily seen to be positive.
\end{proof}

\begin{lemma}\label{A}
For $k>0$ and $x\geq0$, we have
\begin{align*}
\int_{0}^x\,u^{k-1}\,e^{-u} \,du\asymp\left(\frac{x}{1+x}\right)^{k}.
\end{align*}
\end{lemma}

\begin{proof}
The result is clearly true if $0\leq x<1$ (use $e^{-1}\le  e^{-x}\le 1$  and integrate).  If $x\geq1$ then
\begin{align*}
\int_{0}^1\,u^{k-1}\,e^{-u} \,du\leq \int_{0}^x\,u^{k-1}\,e^{-u} \,du<\int_{0}^\infty\,u^{k-1}\,e^{-u} \,du
\end{align*}
and the result follows.
\end{proof}

In a way, the next result contains the essence of the proof of Theorem \ref{trig} for the root system $A_1$.

\begin{lemma}\label{B}
Suppose $a\geq0$.  For $x\geq 0$, we have 
\begin{align*}
\int_0^x\,e^{-a\,u}\,\left(\frac{u}{1+u}\right)^{k-1}\,du\asymp
\left(\frac{x}{1+x}\right)^k \,\frac{1+x}{(1+a\,x)^{k}}
\,\left(\frac{1+a\,(1+x)}{1+a}\right)^{k-1}.
\end{align*}
\end{lemma}

\begin{proof}
Let $A$ represent the left hand side.  
If $0\leq x\leq2$ then by Lemma \ref{A}, we have
\begin{align*}
A\asymp \int_0^x\,e^{-a\,u}\,u^{k-1}\,du
\asymp \left(\frac{x}{1+a\,x}\right)^k
\end{align*}
which gives the result in that case.

If $x\geq 2$ then using Lemma \ref{A} once more and the bound $1-e^{-u}\asymp u/(1+u)$ we have, 
\begin{align*}
A
&\asymp \int_0^1\,e^{-a\,u}\,u^{k-1}\,du+\int_1^x\,e^{-a\,u}\,du
\asymp \left(\frac{1}{1+a}\right)^k+\frac{e^{-a}-e^{-a\,x}}{a}\\
&=\left(\frac{1}{1+a}\right)^k+e^{-a}\,\frac{1-e^{-a\,(x-1)}}{a}
\asymp\left(\frac{1}{1+a}\right)^k+e^{-a}\,\frac{x-1}{1+a\,(x-1)}
\asymp\left(\frac{1}{1+a}\right)^k+e^{-a}\,\frac{x}{1+a\,x}
\end{align*}
which gives the result in that case (consider separately $0\leq a\leq 1$ and $a>1$).
\end{proof}

%%%%%%%%%%%%%%%%%%%%%%%%%%%%%%%%%
\section{Case $A_n$. Proof of Theorem \ref{trig}}\label{P}

\subsection{A recursive formula for spherical functions of type $A_n$}

The following result is an important tool of the proof of Theorem \ref{trig} (see \sout{for example} \cite{Sawyer0}).
\begin{theorem}\label{old}
For $X\in \a^+\subset \R^{n+1}$ and $\lambda\in\R^{n+1}$, we define $\phi_\lambda(e^X)=e^{\lambda(X)}$ when $n=0$ and, for $n\geq 1$,
\begin{align}
\phi_\lambda(e^X)
=\frac{\Gamma(k\,(n+1))}{\Gamma(k)^{n+1}}
e^{\lambda_{n+1}\,\sum_{j=1}^{n+1}\,x_j}
\,\int_{E(X)}\, \phi_{\lambda_0}(e^Y)\,S^{(k)}(Y,X)\,d(Y)^{2\,k}\,dY
\label{Spherical0}
\end{align}
where $E(X)
=\{Y=\diag[y_1,\dots,y_n]\colon x_{j+1}\leq y_j\leq x_j\}$, $\lambda(X)=\sum_{j=1}^{n+1}\,\lambda_j\,x_j$,
$\lambda_0(Y)=\sum_{i=1}^n\,(\lambda_i-\lambda_{n+1})\,y_i$, $d(X)=\prod_{r<s}\,\sinh(x_r-x_s)$, $d(Y)=\prod_{r<s}\,\sinh(y_r-y_s)$ and
\begin{align*}
S^{(k)}(Y,X)&=d(X)^{1-2\,k}\,d(Y)^{1-2\,k}\,\left[\prod_{r=1}^n
\,\left(\prod_{s=1}^r\,\sinh(x_s-y_r)
\,\prod_{s=r+1}^{n+1}\,\sinh(y_r-x_s)\right)\right]^{k-1}.
\end{align*}
Then $\phi_\lambda$ is the (Weyl-invariant) hypergeometric function for the root system $A_n$.
\end{theorem}

\subsection{An equivalent form of Theorem \ref{trig}}

\begin{note}
Define 
\begin{align*}
T_n^{(r)}(X,Y)&=\left( \prod_{s=1}^r\,\frac{x_s-y_r}{1+x_s-y_r}
\,\prod_{s=r+1}^{n+1}\,
\frac{y_r-x_s}{1+y_r-x_s}\right)^{k-1},~1\leq r\leq n\\
T_n(X,Y)&=\prod_{r=1}^{n}
\,T_n^{(r)}(X,Y),\qquad X\in \R^{n+1}, Y\in{\R^n},\\
P_n(\Lambda,Y)&=\prod_{\alpha\in \Sigma_{n-1}^+ }\,\frac{\alpha(Y)\,(1+\alpha(\Lambda)(1+\alpha(Y)))^{k-1}}{(1+\alpha(\Lambda)\alpha(Y))^{k}},
\qquad \Lambda,Y\in \R^n; \qquad P_1=1.
\end{align*}
One can see  $T_n$ as a product of terms of an $n\times (n+1)$ table of factors $\left(\frac{x_s-y_r}{1+x_s-y_r}\right)^{k-1}$ when $s\le r$ and $\left(\frac{y_r-x_s}{1+y_r-x_s}\right)^{k-1}$  when $s> r$ which are all positive when $Y\in E(X).$ Then $T_n^{(r)}(X,Y)$ is the product of the $r$-th row of the table.

\end{note}

\begin{proposition}\label{induc}
Theorem \ref{trig} is equivalent to
\begin{align}\label{eq:I(n)}
I^{(n)}\asymp\frac{\pi(X)^{2\,k-1}}{\prod_{\alpha\in \Sigma_n^+}\,(1+\alpha(\lambda)\alpha(X))^{k}}
\frac{\prod_{\alpha\in \Sigma_n^+}\,(1+\alpha(\lambda)(1+\alpha(X)))^{k-1}}{ \prod_{\alpha\in \Sigma_n^+}\,(1+\alpha(X))^{2k-2}
\prod_{i=1}^n\,(1+\lambda_i-\lambda_{n+1})^{k-1}}, \qquad \lambda,X\in \overline{\a^+},
\end{align}
where, for $\lambda,X\in \R^{n+1}$ 
\begin{align}\label{I(n)}
I^{(n)}(\lambda,X)=
\int_{x_{n+1}}^{x_n}\,\dots\int_{x_2}^{x_1}\,e^{-\sum_{i=1}^n\,(\lambda_i-\lambda_{n+1})\,(x_i-y_i)}\,P_n(\lambda_{1,\ldots,n},Y)\, T_n(X,Y) \,dy_1\dots dy_n.
\end{align}

\end{proposition}

\begin{proof}
This follows using induction and from the fact that for $x\geq 0$, $\sinh x\asymp e^x\,x/(1+x)$ and some simplifications.
\end{proof}

\begin{remark}
If we assume that $\gamma=x_m-x_{m+1}$ is the largest positive root in $X$ then we have either $y_i-y_j\asymp \gamma$ or $y_i-y_j\lesssim \gamma$ for $i<j$ 
(similarly for $x_i-y_j$, $i\leq j$ and $y_i-x_j$, $i<j$). The proof of
the estimate \eqref{eq:I(n)} will be done with the largest positive root $\gamma$
fixed. Moreover, the following result will greatly simplify the proof of the  estimate \eqref{eq:I(n)}.
\end{remark}

\begin{proposition}\label{truncated}
Assume that $\gamma=x_n-x_{n+1}$ is the largest positive root in $X$ and let
\begin{align*}
I_1&=
\int_{M_n}^{x_n}\,\dots\int_{x_2}^{x_1}\,e^{-\sum_{i=1}^n\,(\lambda_i-\lambda_{n+1})\,(x_i-y_i)}\,P_n(\lambda_{1,\ldots,n},Y)\, T_n(X,Y) \,dy_1\dots dy_n
\end{align*}
where $M_n=(x_n+x_{n+1})/2$.  
Then $I_1\asymp I^{(n)}$, when $\lambda$, $X\in \overline{\a^+}$.
\end{proposition}

\begin{proof}
Let $I_2=I^{(n)}-I_1$.
In $I_1$ and $I_2$, consider only the corresponding integral in $y_n$, calling the resulting expressions $\tilde{I}_1$ and  $\tilde{I}_2$. 
To prove the result, it suffices to show that $\tilde{I}_2\lesssim \tilde{I}_1$.

Let $Q_n=(3\,x_n+x_{n+1})/4$.  Observe that for $y_n\in[M_n,Q_n]$, we have
$y_n-x_{n+1}\asymp\gamma$, $x_i-y_n\asymp\gamma$, $1\leq i\leq n$, and $y_i-y_n\asymp\gamma$, $1\leq i\leq n-1$. Thus we have
\begin{align*}
\tilde{I}_1&\gtrsim \int_{M_n}^{Q_n}\,
\,e^{-(\lambda_n-\lambda_{n+1})\,(x_n-y_n)}\,\prod_{i=1}^n
\,\frac{(y_i-y_n)\,(1+(\lambda_i-\lambda_n)\,(1+y_i-y_n))^{k-1}}{(1+(\lambda_i-\lambda_n)\,(y_i-y_n))^{k}}
\\&\qquad
 \left(\prod_{i=1}^n\,\frac{x_i-y_n}{1+x_i-y_n}\,\frac{y_n-x_{n+1}}{1+y_n-x_{n+1}}\right)^{k-1} \,dy_n\\
&\gtrsim e^{-(\lambda_n-\lambda_{n+1})\,\gamma/2}\,\frac{\gamma}4\,\prod_{i=1}^{n-1}\,\frac{\gamma\,(1+(\lambda_i-\lambda_n)\,
(1+\gamma))^{k-1}}{(1+(\lambda_i-\lambda_n)\,\gamma)^{k}}
\, \left(\left(\frac{\gamma}{1+\gamma}\right)^n\,\frac{\gamma}{1+\gamma}\right)^{k-1}.
\end{align*}

On the other hand,  observing that $x_i-y_n\asymp\gamma$, $1\leq i\leq n$, and $y_i-y_n\asymp\gamma$, $1\leq i\leq n-1$, for $y_n\in[x_{n+1}, M_n]$
\begin{align*}
\tilde{I}_2\lesssim e^{-(\lambda_n-\lambda_{n+1})\,\gamma/2}
\,\prod_{i=1}^{n-1}\,\frac{\gamma\,(1+(\lambda_i-\lambda_n)\,(1+\gamma))^{k-1}}{(1+(\lambda_i-\lambda_n)\,\gamma)^{k}}\, \left(\prod_{i=1}^n\,\frac{\gamma}{1+\gamma}\right)^{k-1}
\,\int_{x_{n+1}}^{M_n}\, \left(\frac{y_n-x_{n+1}}{1+y_n-x_{n+1}}\right)^{k-1} \,dy_n.
\end{align*}

Now, for $k\geq1$ using the fact that $u/(1+u)$ is an increasing function, we have
\begin{align*}
\int_{x_{n+1}}^{M_n}\, \left(\frac{y_n-x_{n+1}}{1+y_n-x_{n+1}}\right)^{k-1} \,dy_n\leq\left(\frac{\gamma}{1+\gamma}\right)^{k-1}\, \int_{x_{n+1}}^{M_n} \,dy_n
=\frac{\gamma}{2}\,\left(\frac{\gamma}{1+\gamma}\right)^{k-1}.
\end{align*}

If $0<k\leq1$, we have
\begin{align*}
\int_{x_{n+1}}^{M_n}\, \left(\frac{y_n-x_{n+1}}{1+y_n-x_{n+1}}\right)^{k-1} \,dy_n \le
\left(\frac{1}{1+\gamma}\right)^{k-1} \,\int_{x_{n+1}}^{M_n}\, \left(y_n-x_{n+1}\right)^{k-1}\,dy_n\asymp \left(\frac{1}{1+\gamma}\right)^{k-1} \,\gamma^{k}.
\end{align*}

In both cases, we can conclude that $\tilde{I}_2\lesssim \tilde{I}_1$, 
\end{proof}

%%%%%%%%%%%%%%%%%%%%%%%%%%%%%%%%%%%%%%%

We now prove that Theorem \ref{trig} holds.

\begin{proof}[{Proof of Theorem \ref{trig}}]
We will use Proposition \ref{induc} and use induction on $n$ to show that the estimate of $I^{(n)}$ given in \eqref{eq:I(n)} holds for the root system $A_n$, $n\geq 1$ with root multiplicity $k>0$.

We first prove the result for $n=1$. Let $\alpha=x_1-x_2$. Using Proposition \ref{truncated} and Lemma \ref{B}, we have
\begin{align*}
I^{(1)}&\asymp\int_{M_1}^{x_1}\,e^{-(\lambda_1-\lambda_2)\,(x_1-y)} 
\,\left(\frac{x_1-y}{1+x_1-y}\,\frac{y-x_2}{1+y-x_2}\right)^{k-1}\,dy\\
&\asymp\left(\frac{\alpha}{1+\alpha}\right)^{k-1}\,\int_{M_1}^{x_1}\,e^{-(\lambda_1-\lambda_2)\,(x_1-y)} 
\,\left(\frac{x_1-y}{1+x_1-y}\right)^{k-1}\,dy\\
&=\left(\frac{\alpha}{1+\alpha}\right)^{k-1}\,\int_0^{\alpha/2}\,e^{-(\lambda_1-\lambda_2)\,u} 
\,\left(\frac{u}{1+u}\right)^{k-1}\,du\\
&\asymp \left(\frac{\alpha}{1+\alpha}\right)^{2\,k-1} \,\frac{1+\alpha}{(1+(\lambda_1-\lambda_2)\,\alpha)^{k}}
\,\left(\frac{1+(\lambda_1-\lambda_2)\,(1+\alpha)}{1+(\lambda_1-\lambda_2)}\right)^{k-1}
\end{align*}
which is the desired result by Proposition \ref{induc}.

Assume that the result holds for the root systems $A_1$, $A_2$, \ldots, $A_{n-1}$. 

Fix $1\leq m<n$ and suppose that $\gamma(X)=x_m-x_{m+1}$ is the largest simple positive root  in $X$. We will discuss the case $m=n$ at the end. 
We  divide the integral $I^{(n)}$ in two parts $I_1$ and $I_2$ corresponding to integration in $y_m$ on
the segment $[M_m,x_m]$  and $[x_{m+1}, M_m]$ (recall that $M_m=(x_m+x_{m+1})/2$ and $Q_m=(3\,x_m+x_{m+1})/4$), respectively. The proof consists in two steps:

\bigskip

\noindent {\bf Step 1:}   Show that $I_1$ has the  asymptotics given in \eqref{eq:I(n)}.\\

\noindent {\bf Step 2:}  Show that $I_2\lesssim I_1$.\\

\noindent {\bf   Proof of  Step 1.} Note that for $y_m\in[M_m,x_m]$,  we have 
$x_i-y_j\asymp \gamma$, $i\leq m$, $m<j\leq n$, $y_i-x_j\asymp \gamma$, $i\leq m$, $j\geq m+2$, $y_i-y_j\asymp \gamma$, $i\leq m$, $m<j\leq n$.  It follows that 
\begin{align*}
I_1&=
\int_{x_{n+1}}^{x_n}\dots\int_{M_m}^{x_m}\dots\int_{x_2}^{x_1}\,e^{-\sum_{i=1}^n\,(\lambda_i-\lambda_{n+1})\,(x_i-y_i)}\, T_n(X,Y) \, P_n(\lambda,Y)  
\,dy_1\dots dy_n\\
&=\int_{x_{n+1}}^{x_n}\dots\int_{M_m}^{x_m}\dots\int_{x_2}^{x_1}\,e^{-\sum_{i=1}^n\,(\lambda_i-\lambda_{n+1})\,(x_i-y_i)}\\
&\qquad\, T_m(X_{1,\ldots,m+1}, Y_{1,\ldots,m}) \,T_{n-m}(X_{m+1,\ldots,n+1}, Y_{m+1,\ldots,n}) \, R_1(X,Y)\\
&\qquad\, P_m(\lambda_{1,\ldots,m}, Y_{1,\ldots,m})P_{n-m}(\lambda_{m+1,\ldots,n}, Y_{m+1,\ldots,n}) \,R_2(\lambda,Y)\, dy_1\ldots dy_n,
\end{align*}
and the terms $R_1=T_n/(T_m T_{n-m})$ and
$R_2=P_n/(P_m P_{n-m}) $ have the estimates
\begin{align*}
R_1(X,Y) &\asymp \left(\frac{\gamma}{1+ \gamma}
\right) ^{2m(n-m)(k-1)}=:r_1(X),\\
R_2(\lambda,Y) &\asymp\prod_{i\le m<j\le n}\frac{\gamma((1+(\lambda_i-\lambda_j)(1+\gamma))^{k-1}}{
1+(\lambda_i-\lambda_j)\gamma)^{k}}:=r_2(\lambda,X).
\end{align*}

After replacing the terms $R_1$ and $R_2$ by the estimates $r_1$ and $r_2$, the remaining integrand factorizes
and  by Fubini theorem  and Proposition \ref{truncated}, we get
\begin{align*}
I_1\asymp r_1(X) \,r_2(\lambda,X)
I^{(m)}(\lambda_{1,\ldots,m,n+1}, X_{1,\ldots,m+1})\,
I^{(n-m)}(\lambda_{m+1,\ldots,n+1}, X_{m+1,\ldots,n+1}).
\end{align*}

By the induction hypothesis on $A_m$ and on $A_{n-m}$, we finally obtain
\begin{align*}
I_1&\asymp
\left(\frac{\gamma}{1+ \gamma}
\right) ^{2m(n-m)(k-1)} \,
\prod_{i\le m<j\le n}
\,\frac{\gamma\,((1+(\lambda_i-\lambda_j)(1+\gamma))^{k-1}}{
1+(\lambda_i-\lambda_j)\gamma)^{k}}\\
\, &
\frac{\pi(X_{1,\ldots,m+1})^{2\,k-1}}{\prod_{i<j\le m}
\,(1+(\lambda_i-\lambda_j)(x_i-x_j))^{k}}
\frac{
\prod_{i<j\le m}
\,(1+(\lambda_i-\lambda_j)(1+x_i-x_j))^{k-1}
}{ \prod_{i<j\le m+1}\,(1+x_i-x_j)^{2k-2}
}\\
&\prod_{i=1}^m
\frac{(1+(\lambda_i-\lambda_{n+1})(1+x_i-x_{m+1}))^{k-1}}{(1+(\lambda_i-\lambda_{n+1})(x_i-x_{m+1}))^{k}(1+\lambda_i-\lambda_{n+1})^{k-1} }\\
&
\frac{\pi(X_{m+1,\ldots,n+1})^{2\,k-1}}{\prod_{m+1\le i<j\le n+1}\,(1+(\lambda_i-\lambda_j)(x_i-x_j))^{k}}
\frac{
\prod_{m+1\le i<j\le n+1}
\,(1+(\lambda_i-\lambda_j)(1+x_i-x_j))^{k-1}
}{ \prod_{m+1\le i<j\le n+1}\,(1+x_i-x_j)^{2k-2}
\prod_{i=m+1}^n\,(1+\lambda_i-\lambda_{n+1})^{k-1}}
\end{align*}

Using the fact that $x_i-x_j\asymp \gamma$
when $i\le m$ and $j\ge m+1$, we see that the last expression has the desired asymptotics  \eqref{eq:I(n)}.

\noindent {\bf Proof of Step 2.}  We now show that $I_2=I^{(n)}-I_1\lesssim I_1$.  As before, we show instead that $\tilde{I}_2\lesssim \tilde{I}_1$ 
where $\tilde{I}_1$ (resp.{} $\tilde{I}_2$) represents the portion of $I_1$ (resp.{} $I_2$) where $y_m$ appears. 

Note that for $y_m\in[M_m,Q_m]$,  we have 
$x_i-y_m\asymp \gamma$, $i\leq m$ and $y_m-x_j\asymp \gamma$, $j>m$, and $|y_i-y_m|\asymp \gamma$, $i\not=m $. 

It follows that 
\begin{align*}
\tilde{I}_1& \gtrsim
\int_{M_m}^{Q_m} \,e^{-(\lambda_m-\lambda_{n+1})\,(x_m-y_m)} \, T_n^{(m)}(X,Y) 
\, \prod_{i=1}^{m-1}\,\frac{(y_i-y_m)\,(1+(\lambda_i-\lambda_m)\,(1+y_i-y_m))^{k-1}}{(1+(\lambda_i-\lambda_m)\,(y_i-y_m))^{k}} 
\\&\qquad
\prod_{i=m+1}^n\,\frac{(y_i-y_m)\,(1+(\lambda_m-\lambda_i)\,(1+y_m-y_i))^{k-1}}{(1+(\lambda_m-\lambda_i)\,(y_m-y_i))^{k}} 
\,dy_m \\
&\gtrsim 
e^{-(\lambda_m-\lambda_{n+1})\,\gamma/2}\,\frac{\gamma}{4}
\, 
\left(\frac{\gamma}{1+\gamma}\right)^{(n+1)\,(k-1)}\,
\prod_{i=1}^{m-1}\,\frac{\gamma\,(1+(\lambda_i-\lambda_m)\,(1+\gamma))^{k-1}}{(1+(\lambda_i-\lambda_m)\,\gamma)^{k}} 
\\&\qquad
\,\prod_{i=m+1}^n\,\frac{\gamma\,(1+(\lambda_m-\lambda_i)\,(1+\gamma))^{k-1}}{(1+(\lambda_m-\lambda_i)\,\gamma)^{k}}.
\end{align*}

On the other hand, since $x_i-y_m\asymp \gamma$, $i\leq m$, and $y_i-y_m\asymp \gamma$, $i<m$, when $y_m\in[x_{m+1},M_m]$, we have
\begin{align}
\tilde{I}_2& \lesssim
\,e^{-(\lambda_m-\lambda_{n+1})\,\gamma/2} \,\int_{x_{m+1}}^{M_m}\, T_n^{(m)}(X,Y) 
\, \prod_{i=1}^{m-1}\,\frac{(y_i-y_m)\,(1+(\lambda_i-\lambda_m)\,(1+y_i-y_m))^{k-1}}{(1+(\lambda_i-\lambda_m)\,(y_i-y_m))^{k}} 
\nonumber\\&\qquad
\prod_{i=m+1}^n\,\frac{(y_m-y_i)\,(1+(\lambda_m-\lambda_i)\,(1+y_m-y_i))^{k-1}}{(1+(\lambda_m-\lambda_i)\,(y_m-y_i))^{k}} 
\,dy_m \nonumber\\\
& \lesssim
\,e^{-(\lambda_m-\lambda_{n+1})\,\gamma/2} 
\,\left(\frac{\gamma}{1+\gamma}\right)^{m\,(k-1)}\, \prod_{i=1}^{m-1}\,\frac{\gamma\,(1+(\lambda_i-\lambda_m)
\,(1+\gamma))^{k-1}}{(1+(\lambda_i-\lambda_m)\,\gamma)^{k}} 
\label{it}
\\&
\,\int_{x_{m+1}}^{M_m}\,\prod_{i=m+1}^{n+1}\,\left(\frac{y_m-x_i}{1+y_m-x_i}\right)^{k-1}
\,\prod_{i=m+1}^n\,\frac{(y_m-y_i)\,(1+(\lambda_m-\lambda_i)\,(1+y_m-y_i))^{k-1}}{(1+(\lambda_m-\lambda_i)\,(y_m-y_i))^{k}} 
\,dy_m.\nonumber
\end{align}

If $k\leq 1$, referring to the function $F_2$ of Lemma \ref{incr}, we have from \eqref{it}
\begin{align*}
\tilde{I}_2& \lesssim
e^{-(\lambda_m-\lambda_{n+1})\,\gamma/2} 
\left(\frac{\gamma}{1+\gamma}\right)^{m\,(k-1)}
\, \prod_{i=1}^{m-1}\,\frac{\gamma\,(1+(\lambda_i-\lambda_m)
\,(1+\gamma))^{k-1}}{(1+(\lambda_i-\lambda_m)\,\gamma)^{k}} 
\,
\\&\prod_{i=m+1}^{n+1}\,\left(\frac{1}{1+\gamma}\right)^{k-1}
\,\prod_{i=m+1}^n\,\frac{\gamma\,(1+(\lambda_m-\lambda_i)\,(1+\gamma))^{k-1}}{(1+(\lambda_m-\lambda_i)\,\gamma)^{k}} 
\,\int_{x_{m+1}}^{M_m}\,\prod_{i=m+1}^{n+1}\,(y_m-x_i)^{k-1}\,dy_m\\
&\lesssim e^{-(\lambda_m-\lambda_{n+1})\,\gamma/2} 
\left(\frac{\gamma}{1+\gamma}\right)^{m\,(k-1)}
\, \prod_{i=1}^{m-1}\,\frac{\gamma\,(1+(\lambda_i-\lambda_m)\,
(1+\gamma))^{k-1}}{(1+(\lambda_i-\lambda_m)\,\gamma)^{k}} 
\,
\\&\prod_{i=m+1}^{n+1}\,\left(\frac{1}{1+\gamma}\right)^{k-1}
\,\prod_{i=m+1}^n\,\frac{\gamma\,(1+(\lambda_m-\lambda_i)\,(1+\gamma))^{k-1}}{(1+(\lambda_m-\lambda_i)\,\gamma)^{k}} 
\,\int_{x_{m+1}}^{M_m}\,(y_m-x_{m+1})^{(n+1-m)\,(k-1)}\,dy_m\\
&\asymp  e^{-(\lambda_m-\lambda_{n+1})\,\gamma/2} 
\,\left(\frac{\gamma}{1+\gamma}\right)^{m\,(k-1)}
\, \prod_{i=1}^{m-1}\,\frac{\gamma\,(1+(\lambda_i-\lambda_m)
\,(1+\gamma))^{k-1}}{(1+(\lambda_i-\lambda_m)\,\gamma)^{k}} 
\\&\prod_{i=m+1}^{n+1}\,\left(\frac{1}{1+\gamma}\right)^{k-1}
\,\prod_{i=m+1}^n\,\frac{\gamma\,(1+(\lambda_m-\lambda_i)\,(1+\gamma))^{k-1}}{(1+(\lambda_m-\lambda_i)\,\gamma)^{k}} 
\,\gamma^{(n+1-m)\,(k-1)+1}\lesssim \tilde{I}_1.
\end{align*}

If $k\geq 1$, using $y_m-x_i\leq y_m-y_i$, referring to the functions $F_1$ and $F_3$ of Lemma \ref{incr} and rewriting \eqref{it}, we have
\begin{align*}
\tilde{I}_2& \lesssim
\,e^{-(\lambda_m-\lambda_{n+1})\,\gamma/2} 
\,\left(\frac{\gamma}{1+\gamma}\right)^{m\,(k-1)}
\, \prod_{i=1}^{m-1}\,\frac{\gamma\,(1+(\lambda_i-\lambda_m)\,(1+\gamma))^{k-1}}{(1+(\lambda_i-\lambda_m)\,\gamma)^{k}} 
\\&
\,\int_{x_{m+1}}^{M_m}\,
\left(\frac{y_m-x_{n+1}}{1+y_m-x_{n+1}}
\right)^{k-1}
\,\prod_{i=m+1}^n\,
\frac{y_m-y_i}{1+(\lambda_m-\lambda_i)\,(y_m-y_i)}
\\&
\,\prod_{i=m+1}^n\,\left(\frac{(y_m-y_i)\,(1+(\lambda_m-\lambda_i)\,(1+y_m-y_i))}{(1+y_m-y_i)\,(1+(\lambda_m-\lambda_i)\,(y_m-y_i))} \right)^{k-1}
\,dy_m\\
& \lesssim
\,e^{-(\lambda_m-\lambda_{n+1})\,\gamma/2} 
\,\left(\frac{\gamma}{1+\gamma}\right)^{m\,(k-1)}
\, \prod_{i=1}^{m-1}\,\frac{\gamma\,(1+(\lambda_i-\lambda_m)
\,(1+\gamma))^{k-1}}{(1+(\lambda_i-\lambda_m)\,\gamma)^{k}} 
\\&
\,\frac{\gamma}2\,
\left(\frac{\gamma}{1+\gamma}
\right)^{k-1}
\,\prod_{i=m+1}^n\,
\frac{\gamma}{1+(\lambda_m-\lambda_i)\,\gamma}
\,\prod_{i=m+1}^n\,\left(\frac{\gamma\,(1+(\lambda_m-\lambda_i)\,(1+\gamma))}{(1+\gamma)\,(1+(\lambda_m-\lambda_i)\,\gamma)} \right)^{k-1}\\
&=\,e^{-(\lambda_m-\lambda_{n+1})\,\gamma/2} 
\,\left(\frac{\gamma}{1+\gamma}\right)^{(n+1)\,(k-1)}
\, \prod_{i=1}^{m-1}\,\frac{\gamma\,(1+(\lambda_i-\lambda_m)\,(1+\gamma))^{k-1}}{(1+(\lambda_i-\lambda_m)\,\gamma)^{k}} 
\\&
\,\frac{\gamma}2
\,\prod_{i=m+1}^n\,\frac{\gamma\,(1+(\lambda_m-\lambda_i)
\,(1+\gamma))^{k-1}}{(1+(\lambda_m-\lambda_i)\,\gamma)^k} 
\lesssim \tilde{I}_1.
\end{align*}

To conclude, we reason by symmetry, as explained below. By the structure of the root system $A_n$, the case $\alpha_n$ maximal is equivalent to the case $\alpha_1$ maximal.  Indeed, in formula 
\eqref{Spherical0}, one does not assume that $\lambda\in\overline{\a^+}$.  We also know that $\phi_\lambda(e^X)$ is invariant under permutation of its $\lambda$ argument.  Hence one can re-write \eqref{Spherical0}
by exchanging $\lambda_1$ and $\lambda_{n+1}$,
\begin{align*}
\phi_\lambda(e^X)&=e^{\lambda(X)}\ \hbox{if $n=1$ and}\nonumber\\
\phi_\lambda(e^X)
&=\frac{\Gamma(k\,(n+1))}{\Gamma(k)^{n+1}}
e^{\lambda_1\,\sum_{j=1}^{n+1}\,x_j}
\int_{E(X)}\, \phi_{\widetilde{\lambda_0}}(e^Y)\,S^{(k)}(Y,X)\,d(Y)^{2\,k}\,dY
\end{align*}
where $\widetilde{\lambda_0}(Y)=\sum_{r=2}^{n+1}\,(\lambda_r-\lambda_1)\,y_{r-1}$.

We used the fact that
\begin{align*}
\phi_{[\lambda_{n+1}-\lambda_1,\lambda_2-\lambda_1,\dots,\lambda_n-\lambda_1]}(e^Y)
=\phi_{[\lambda_2-\lambda_1,\dots,\lambda_n-\lambda_1,\lambda_{n+1}-\lambda_1]}(e^Y).
\end{align*}

Theorem \ref{trig} is then equivalent to
\begin{align*}
J^{(n)}\asymp\frac{\pi(X)^{2\,k-1}}{\prod_{\alpha\in \Sigma_n^+}\,(1+\alpha(\lambda)\alpha(X))^{k}}
\frac{\prod_{\alpha\in \Sigma_n^+}\,(1+\alpha(\lambda)(1+\alpha(X)))^{k-1}}{ \prod_{\alpha\in \Sigma_n^+}\,(1+\alpha(X))^{2k-2}
\prod_{i=2}^{n+1}\,(1+\lambda_1-\lambda_i)^{k-1}}
\end{align*}
where, for $\lambda,X\in \R^{n+1}$,
\begin{align*}
J^{(n)}(\lambda,X)=
\int_{x_{n+1}}^{x_n}\,\dots\int_{x_2}^{x_1}\,e^{-\sum_{i=1}^n\,(\lambda_1-\lambda_{i+1})\,(y_i-x_{i+1})}\,P_n(\lambda_{2,\ldots,n+1},Y)\, T_n(X,Y) \,dy_1\dots dy_n.
\end{align*}
The term $J^{(n)}$ corresponds to a constant multiple of $e^{-\lambda(X)}\,d(X)^{2\,k-1}\,\phi_\lambda(e^X)$ in which we have replaced $\phi_{\widetilde{\lambda_0}}(e^Y)$ by its asymptotic expression proposed in Theorem \ref{trig}.  One then proves the case $\alpha_n$ maximal as one proves the case $\alpha_1$ maximal.

This concludes the proof of the Theorem \ref{trig} for $X\in\a^+$ (recall that the formula \eqref{Spherical0} holds for $X\in\a^+$). The estimates that we find for $\phi_\lambda(e^X)$ extend to $X\in\overline{\a^+}$ by continuity.
\end{proof}

\section{{Other evidence for Conjecture \ref{conj:all}}}\label{evidence}

\subsection{Comparison with known estimates} \label{comp}

Recall the estimates from Narayanan and al in \cite{Narayanana} (refer also to \cite{Schapira}):
\begin{align}
C_1(\lambda)\,e^{\lambda-\rho}(X)\,\prod_{\alpha\in\Sigma^{++}_{\lambda}}\,(1+\alpha(X)) \leq \phi_\lambda(X)\leq C_2(\lambda)\,e^{\lambda-\rho}(X)\,\prod_{\alpha\in\Sigma^{++}_{\lambda}}\,(1+\alpha(X))
\label{N}
\end{align}
where, as before, 
$\Sigma^{++}_{\lambda}=\{\alpha\in  \Sigma^{++}\colon \alpha(\lambda)=0\}$.

We will show that our conjecture is consistent with the bound \eqref{N}.  We will need a technical Lemma.

\begin{lemma}\label{ab}
Assume $u\geq 0$ and $a>0$.  Then
\begin{align*}
\frac{1}{1+a}\leq \frac{1+u}{(1+a\,u)^{k}} \,(1+(1+u)\,a)^{k-1}\leq \frac{(1+a)^k}{a}.
\end{align*}
\end{lemma}

\begin{proof}
Refer to Lemma \ref{incr}.  We have
\begin{align*}
f(u)&=\frac{1+u}{(1+a\,u)^{k}} \,(1+(1+u)\,a)^{k-1}
=\frac{1+u}{1+(1+u)\,a} \,\left(\frac{1+(1+u)\,a}{1+a\,u}\right)^k=F_5(u)/F_4(u)^k.
\end{align*}

The function $F_5(u)$ increases in $u$ and therefore $1/(1+a)=F_5(0)\leq F_5(u)\leq F_5(\infty)=1/a$.

The function  $1/F_4(u)$ decreases in $u$ and therefore $1+a=1/F_4(0)\geq 1/F_4(u)\geq 1/F_4(\infty)=1$.
\end{proof}

\begin{proposition}\label{Nara}
Conjecture \ref{conj:all} is consistent with \eqref{N}.
\end{proposition}

\begin{proof}
Assume the bound proposed in Conjecture \ref{conj:all}:  
\begin{align*}
\phi_\lambda(e^X)&\asymp 
e^{(\lambda-\rho)(X)}\,\prod_{\alpha\in\Sigma^{++}}\,\frac{1+\alpha(X)}{(1+\alpha(\lambda)\,\alpha(X))^{k(\alpha)+k(2\,\alpha)}}
\,\left(\frac{1+\alpha(\lambda)\,(1+\alpha(X))}{1+\alpha(\lambda)}\right)^{k(\alpha)+k(2\,\alpha)-1}\\
&=e^{(\lambda-\rho)(X)}\,\prod_{\alpha\in\Sigma^{++},\alpha(\lambda)=0}\,\frac{1+\alpha(X)}{(1+\alpha(\lambda)\,\alpha(X))^{k(\alpha)+k(2\,\alpha)}}
\,\left(\frac{1+\alpha(\lambda)\,(1+\alpha(X))}{1+\alpha(\lambda)}\right)^{k(\alpha)+k(2\,\alpha)-1}
\\&\qquad
\,\prod_{\alpha\in\Sigma^{++},\alpha(\lambda)>0}\,\frac{1+\alpha(X)}{(1+\alpha(\lambda)\,\alpha(X))^{k(\alpha)+k(2\,\alpha)}}
\,\left(\frac{1+\alpha(\lambda)\,(1+\alpha(X))}{1+\alpha(\lambda)}\right)^{k(\alpha)+k(2\,\alpha)-1}\\
&=e^{(\lambda-\rho)(X)}\,\prod_{\alpha\in\Sigma^{++}},\alpha(\lambda)=0\,(1+\alpha(X))
\\&\qquad
\,\prod_{\alpha\in\Sigma^{++},\alpha(\lambda)>0}\,\frac{1+\alpha(X)}{(1+\alpha(\lambda)\,\alpha(X))^{k(\alpha)+k(2\,\alpha)}}
\,\left(\frac{1+\alpha(\lambda)\,(1+\alpha(X))}{1+\alpha(\lambda)}\right)^{k(\alpha)+k(2\,\alpha)-1}.
\end{align*}

In order to show that this is consistent with \eqref{N}, we only have to show that each term 
\begin{align*}
\frac{1+\alpha(X)}{(1+\alpha(\lambda)\,\alpha(X))^{k(\alpha)+k(2\,\alpha)}}
\,(1+\alpha(\lambda)\,(1+\alpha(X)))^{k(\alpha)+k(2\,\alpha)-1}
\end{align*}
is bounded below and above by expressions only depending on $\lambda$ whenever $\alpha(\lambda)>0$ for a positive root $\alpha $.  This follows from Lemma \ref{ab}.

\end{proof}

\subsection{$BC_1$ case}
Recall that the only rank 1 root systems are $A_1$ and $BC_1$.  In this section we discuss Conjecture \ref{conj:all} for the system $BC_1$. Denote $k_1=k(\alpha)$ and $k_2=k(2\alpha)$. We have $\rho=k_1+2k_2$.  

Conjecture \ref{conj:all} reads in this case
\begin{align*}
\phi_\lambda(e^t)
\asymp 
e^{(\lambda-\rho)(t)}\,\frac{1+t}{(1+\lambda\,t)^{k_1+k_2 }}
\,\left(\frac{1+\lambda\,(1+t)}{1+\lambda}\right)^{
k_1+k_2-1}
\end{align*}

We need to prove four following  asymptotics, with the notation $k=k_1+k_2$:
\begin{align}
e^{-(\lambda-\rho)(t)}\,\phi_\lambda(e^t) &\asymp\left\lbrace
\begin{array}{cl}
1+t&\hbox{if $\lambda t\leq1$},\\
 (\lambda\,t)^{-k} &\hbox{if $\lambda t\geq1$,~$t\leq 1$},\\
\lambda^{-1}&\hbox{if $\lambda t\geq1$, $t\geq 1$, $\lambda\leq1$},\\
 \lambda^{-k}&\hbox{if $\lambda t\geq1$, $ t\geq 1$, $\lambda\geq1$}.
\end{array}
\right.\label{Four}
\end{align}

\begin{lemma}\label{log}
For $0\leq r\leq 1$ and $0\leq t\leq 1$, there exists $C>0$ independent of $t$ and $r$ such that
\begin{align*}
\left|\log(\cosh t+r\,\sinh t)-r\,t\right|\leq C\,t^2.
\end{align*}
\end{lemma}

\begin{proof}
Let $F(r)=\log(\cosh t+r\,\sinh t) -r\,t$.  We find that $F'(r)=0$ only if $r=r_0=(\sinh t-t\,\cosh t)/(t\,\sinh t)$.  The maximum and minimum of $F(r)$ on $[0,1]$ can only occur 
if $r=0$, $r=1$ or $r=r_0$.  This corresponds to the values of $F(r)$ equal to
$\log(\cosh(t))/t^2$, 0 or $(\cosh(t)\,t + \log(\sinh t/t)\,\sinh t - \sinh t)/(t^2\,\sinh t)$.  The result follows.

\end{proof}

\begin{lemma}\label{special}
Assume $x_1\geq x_2$.  Let $\tilde{\phi}^{(k)}_\lambda$ denote the spherical function for the  $A_1$ root system.  Then
\begin{align*}
\int_{x_2}^{x_1}\,e^{\mu\,y}\,((e^{2\,x_1}-e^{2\,y})\,(e^{2\,y}-e^{2\,x_2}))^{k-1}\,dy
&=C_k\,e^{(k-1)\,(x_1+x_2)}\,\sinh^{2\,k-1}(x_1-x_2)\,\tilde{\phi}^{(k)}_{[\mu+2\,(k-1),0]}(e^X).
\end{align*}
\end{lemma}

\begin{proof}
Note that $(e^{2\,x_1}-e^{2\,y})\,(e^{2\,y}-e^{2\,x_2})=4\,e^{x_1+x_2}\,e^{2\,y}\,\sinh(x_1-y)\,\sinh(y-x_2)$.  
\end{proof}

\begin{proposition}
The spherical functions in the $BC_1$ case satisfy Conjecture \ref{conj:all}.
\end{proposition}

\begin{proof}
The spherical functions for the $BC_1$ case are given by \cite[(5.28)]{Koornwinder} (where we corrected a  small misprint):
\begin{align}
\phi_\lambda(a_t)&=\frac{2\,\Gamma(k_1+k_2+1/2)}{\Gamma(1/2)\,\Gamma(k_1)\,\Gamma(k_2)}
\,\int_0^1\,\int_0^\pi\,|\cosh(t)+r\,e^{i\,\phi}\,\sinh(t)|^{\lambda-\rho}
\,(1-r^2)^{k_1-1}\,r^{2\,k_2-1}\,\sin^{2\,k_2-1}\phi\,r\,dr\,d\phi\label{Koor}
\end{align}
with $\rho=k_1+2\,k_2$.

We divide the region $(\lambda,X)\in\overline{\a^+}\times \overline{\a^+}$ in Regions $I$, $I$, $II$ and $IV$ based on \eqref{Four} (with some technical variations).  Figure \ref{R} in the Introduction illustrates these regions.

\noindent{\bf Region I:} Suppose $0\leq \lambda\,t\leq 1$. Since, for $0\le r \le 1$,
\begin{align*}
e^{-t}= \cosh(t)-\sinh(t)\leq |\cosh(t)+r\,e^{i\,\phi}\,\sinh(t)|\leq \cosh(t)+\sinh(t)=e^t,
\end{align*}
\begin{align*}
\phi_\lambda(a_t)\asymp \phi_0(a_t)\asymp e^{-\rho\,t}\,(1+t)\asymp e^{(\lambda-\rho)\,t}\,(1+t)
\end{align*}
(refer to \eqref{N}) which proves the proposition in this case.

Now,
\begin{align*}
|\cosh(t)+r\,e^{i\,\phi}\,\sinh(t)|^2&=[\cosh(t)+r\,e^{i\,\phi}\,\sinh(t)]\,[\cosh(t)+r\,e^{-i\,\phi}\,\sinh(t)]\\
&=\cosh^2t+2\,r\,\cos\phi \,\sinh t\,\cosh t+r^2\,\sinh^2t\\
&=\cosh^2t+r\,\cos\phi \,\sinh (2\,t)+r^2\,\sinh^2t.
\end{align*}

Hence, using $e^{2\,x}=\cosh^2t+r\,\cos\phi 
\,\sinh (2\,t)+r^2\,\sinh^2t$ and noting that
\begin{align*}
\cos\phi&=\frac{e^{2\,x}-\cosh^2 t-r^2\,\sinh^2t}{r\,\sinh(2\,t)},\\
\sin\phi&=\sqrt{1-\cos^2\phi}=\frac{\left(\left(\cosh t+r\,\sinh t\right)^{2}-e^{2 x}\right)^{1/2} \left(e^{2 x}-\left(\cosh t-r\,\sinh t\right)^{2}\right)^{1/2}}{r\,\sinh(2\,t)},~0\leq \phi\leq \pi,
\end{align*}
we have (the constant $C$ may vary from line to line):

denoting $x_1(r,t)=\log(\cosh t+r\,\sinh t)$ and $x_2(r,t)=\log(\cosh t-r\,\sinh t)$,
\begin{align}
\phi_\lambda(a_t)
&=C\,\int_0^1\,\int_0^\pi\,[\cosh^2t+r\,\cos\phi \,\sinh (2\,t)+r^2\,\sinh^2t]^{(\lambda-\rho)/2}
\,(1-r^2)^{k_1-1}\,r^{2\,k_2-1}\,\sin^{2\,k_2-1}\phi\,r\,dr\,d\phi\nonumber\\
&=\frac{C}{\sinh(2\,t)}\,\int_0^1\,\int_{\log(\cosh t + r\,\sinh t)}^{\log(\cosh t - r\,\sinh t)}\,e^{(\lambda-\rho)\,x}
\,(1-r^2)^{k_1-1}\,r^{2\,k_2-1}\,\sin^{2\,k_2-2}\phi\,\overbrace{(-r\,\sin\phi\sinh  (2\,t))/2\,d\phi}^{e^{2\,x}\,dx}\,dr\label{0}\\
&=\frac{C}{\sinh^{2\,k_2-1}(2\,t)}\,\int_0^1\,\int_{\log(\cosh t - r\,\sinh t)}^{\log(\cosh t + r\,\sinh t)}\,e^{(\lambda+2-\rho)\,x}
\,(1-r^2)^{k_1-1}\,r
\nonumber\\&\qquad\qquad\qquad\qquad\qquad\,
\left[
\left(\left(\cosh t+r\,\sinh t\right)^{2}-e^{2 x}\right) \left(e^{2 x}-\left(\cosh t-r\,\sinh t\right)^{2}\right)
\right]^{k_2-1}\,dr\,dx\nonumber\\
&=\frac{C}{\sinh^{2\,k_2-1}(2\,t)}\,\int_0^1\,(1-r^2)^{k_1-1}
\,e^{(k_2-1)\,(x_1(r,t)+x_2(r,t))}
\nonumber\\&\qquad
\,\sinh^{2\,k_2-1}(x_1(r,t)-x_2(r,t))
\,\tilde{\phi}^{(k_2)}_{[\lambda-\rho+2+2\,(k_2-1),0]}(x_1(r,t),x_2(r,t))
\,r\,dr\nonumber\\
&=\frac{C}{\sinh^{2\,k_2-1}(2\,t)}\,\int_0^1\,(1-r^2)^{k_1-1}
\,(\cosh^2t-r^2\,\sinh^2t)^{k_2-1}
\,\left(\frac{r\,\sinh(2\,t)}{\cosh^2t-r^2\,\sinh^2t}\right)^{2\,k_2-1}
\nonumber\\&\qquad
\,\tilde{\phi}^{(k_2)}_{[\lambda-\rho+2\,k_2,0]}(x_1(r,t),x_2(r,t))
\,r\,dr\nonumber\\
&=C\,\int_0^1\,(1-r^2)^{k_1-1}
\,(\cosh^2t-r^2\,\sinh^2t)^{-k_2}
\,\tilde{\phi}^{(k_2)}_{[\lambda-k_1,0]}(x_1(r,t),x_2(r,t))
\,r^{2\,k_2}\,dr\nonumber
\end{align}
since
\begin{align*}
\sinh(x_1(r,t)-x_2(r,t))=\frac{r\,\sinh(2\,t)}{\cosh^2t-r^2\,\sinh^2t}.
\end{align*}

Remark that
\begin{align}
-t=\log(\cosh t - \sinh t)\leq\log(\cosh t - r\,\sinh t)\leq \log(\cosh t + r\,\sinh t)\leq \log(\cosh t + \sinh t)=t.\label{T}
\end{align}

\noindent{\bf Region II:} Suppose now that $ \lambda\,t\geq 1$ and $0\leq t\leq T_0=\min\{\log 2,1/(2\,k_1),1/(4\,C)\}$ where $C$ is as in Lemma \ref{log}.  In that case,  $\lambda\geq 2\,k_1$,

\begin{align*}
1=\cosh^2 t-\sinh^2t&\leq\cosh^2 t-r^2\,\sinh^2t\leq \cosh^2t\leq \cosh^2(1/2),\\
r\,t&\leq x_1(r,t)-x_2(r,t)\leq 2\,r\,t
\end{align*}
(for the two last inequalities we studied variations of convenient functions and used $1-\sinh^2t\geq 0$ for $t\le \log 2$).

Applying the estimates for $A_1$ with the multiplicity $k_2$ we get
\begin{align*}
\phi_\lambda(a_t)
&\asymp\int_0^1\,(1-r^2)^{k_1-1}\,(\cosh^2t-r^2\,\sinh^2t)^{-k_2}
\,e^{(\lambda-k_1)\,x_1(r,t)}\,
e^{-k_2\,(x_1(r,t)-x_2(r,t))}
\\&\qquad
\,\frac{1+x_1(r,t)-x_2(r,t)}{(1+(\lambda-k_1)\,(x_1(r,t)-x_2(r,t)))^{k_2}}
\\&\qquad
\,\left(\frac{1+(\lambda-k_1)\,(1+x_1(r,t)-x_2(r,t))}{1+\lambda-k_1}\right)^{k_2-1}
\,r^{2\,k_2}\,dr\\
&\asymp\int_0^1\,(1-r)^{k_1-1}\,\,e^{\lambda\,x_1(r,t)}\,\frac{1+r\,t}{(1+\lambda\,(r\,t))^{k_2}}
\,\left(\frac{1+\lambda\,(1+r\,t)}{1+\lambda}\right)^{k_2-1}
\,r^{2\,k_2}\,dr~\text{using \eqref{T}}\\
&\asymp e^{\lambda t}\,\int_0^1\,(1-r)^{k_1-1}\,e^{-\lambda\,[t-\log(\cosh t+r\,\sinh t)]}\,(1+\lambda\,r\,t)^{-k_2}\,r^{2\,k_2}\,dr.
\end{align*}

We have
\begin{align*}
\lefteqn{\int_0^1\,(1-r)^{k_1-1}
\,e^{-\lambda\,[t-\log(\cosh t+r\,\sinh t)]}\,(1+\lambda\,r\,t)^{-k_2}\,r^{2\,k_2}\,dr}\\
&\asymp (1+\lambda\,t)^{-k_2}\,\int_{1/2}^1\,(1-r)^{k_1-1}\,\,e^{-\lambda\,[t-\log(\cosh t+r\,\sinh t)]}\,dr
\\&\qquad
+\int_0^{1/2}\,e^{-\lambda\,[r\,t-\log(\cosh t+r\,\sinh t)+(1-r)\,t]}\,(1+\lambda\,r\,t)^{-k_2}\,r^{2\,k_2}\,dr.
\end{align*}

By Lemma \ref{log},
\begin{align*}
\lefteqn{\int_0^{1/2}\,\,e^{-\lambda\,[r\,t-\log(\cosh t+r\,\sinh t)+(1-r)\,t]}\,(1+\lambda\,r\,t)^{-k_2}\,r^{2\,k_2}\,dr}\\
&\lesssim \int_0^{1/2}\,\,e^{-\lambda\,[(1-r)\,t-C\,t^2]}\,(1+\lambda\,r\,t)^{-k_2}\,r^{2\,k_2}\,dr\\
&\lesssim \int_0^{1/2}\,\,e^{-\lambda\,t\,[1/2-C\,t]}\,(1+\lambda\,r\,t)^{-k_2}\,r^{2\,k_2}\,dr\\
&\lesssim \int_0^{1/2}\,e^{-\lambda\,t/4}\,(1+\lambda\,r\,t)^{-k_2}\,r^{2\,k_2}\,dr
\lesssim (1/2)^{2\,k_2+1}\,e^{-\lambda\,t/4}.
\end{align*}

On the other hand, using $u=t-\log(\cosh t+r\,\sinh t)$,
\begin{align*}
\lefteqn{\int_{1/2}^1\,(1-r)^{k_1-1}\,e^{-\lambda\,[t-\log(\cosh t+r\,\sinh t)]}\,dr}\\
&=\int_0^{t-\log(\cosh t+(\sinh t)/2)}\,e^{-\lambda\,u}\,\left(\frac{2\,(1-e^{-u})}{1-e^{-2\,t}}\right)^{k_1-1}\,\frac{e^{-u}}{1-e^{-2\,t}}\,du\\
&=\frac{C}{(1-e^{-2\,t})^{k_1}}\,\int_0^{t-\log(\cosh t+(\sinh t)/2)}\,e^{-\lambda\,u}\,e^{-u}\,\left(1-e^{-u}\right)^{k_1-1}\,du\\
&\asymp t^{-k_1}\,\int_0^{t-\log(\cosh t+(\sinh t)/2)}\,e^{-\lambda\,u}\,u^{k_1-1}\,du\asymp  t^{-k_1}\,\lambda^{-k_1}\,\int_0^{\lambda\,(t-\log(\cosh t+(\sinh t)/2))}\,e^{-v}\,v^{k_1-1}\,dv\\
&\asymp  t^{-k_1}\,\lambda^{-k_1}
\end{align*}
since
\begin{align*}
\lambda\,(t-\log(\cosh t+(\sinh t)/2)) \geq \lambda\,t/8\geq 1/8
\end{align*}
whenever $0\leq t\leq \log 2$.  This proves the proposition for this case.

\noindent{\bf Region III:}  Suppose that $\lambda t\geq 1$, $t\geq t_0$, $\lambda\leq \min(\frac{\rho}2, \frac12)$.
The constant $t_0$ will be defined in the proof.

According to  \cite[Ex.{} 8, p.{} 484]{Helgason} (see also \cite[p.{} 325]{Opdam} and \cite[p.{} 109]{Shimeno}), 
\begin{align}\label{eq:shimeno}
\phi_\lambda(e^t)= \;_2 F_1
\left( \frac{\rho}2 + \frac{\lambda}2, 
\frac{\rho}2 - \frac{\lambda}2; k_1+k_2
 + \frac12 ; - \sinh^2 t\right)
\end{align}

We apply the formula
\cite[15.8.2]{NIST}, with $z=-\sinh^2 t$ and we get:
\begin{align}
\frac{-\sin(\pi\lambda)}{\pi \Gamma(k_1+k_2+\frac12)}
\phi_\lambda(e^t)
 &=
 \frac{(\sinh t)^{-(\rho+\lambda)}}{\Gamma(\frac{\rho-\lambda}2 )  \Gamma(\frac{k_1+1-\lambda}2)
 \Gamma(\lambda+1)
 }
 \;_2 F_1 (\frac{\rho+\lambda}{2},
  \frac{\lambda-k_1+1}{2}, \lambda+1; \frac1z)\nonumber\\
  &+
 \frac{(\sinh t)^{\lambda-\rho}}{\Gamma(\frac{\rho+\lambda}2 )  \Gamma(\frac{k_1+1+\lambda}2)
 \Gamma(-\lambda+1)
 }
 \;_2 F_1 (\frac{\rho-\lambda}{2},
  \frac{-\lambda-k_1+1}{2}, -\lambda+1; \frac1z) \label{F}
 \end{align}
 
 All the Gamma functions are bounded and bounded away from zero.
 
 As $|\frac1z|<1$, the hypergeometric functions on the right hand side of \eqref{F} are equal to  the hypergeometric power series
\begin{align*}
 {}_2F_1(a,b;c;w)=\sum_{n=0}^\infty  \frac{(a)_n(b)_n}{(c)_n} w^n
 \end{align*}
 with $w=\frac1z$.
Observing that
if $|a|\leq a_0$,  $|b|\leq b_0$ and $|c|\geq c_0>0$ and $w\leq 1/2$ then by manipulating the defining power series, one finds that
\begin{align*}
\left|{}_2F_1(a,b;c;w)-1\right|\leq \frac{a_0\,b_0}{c_0}\,{}_2F_1(a_0+1,b_0+1;c_0+1;1/2)\,w.
\end{align*}
This will allow us to show that the hypergeometric terms on the right hand side of \eqref{F} are bounded and bounded away from 0, if $w=1/z$ is small enough.
Indeed, let us take  
\begin{align*}
w_0=\min\left\{1/2, \frac{1/2}{(\rho/2+1/4)\,(3/2+k_1)/2}\,\frac{1}{{}_2F_1(\rho/2+5/4,k_1/2+7/4, 3/2;1/2)}\right\}.
\end{align*}
In the assumptions of this case, we must suppose $t\ge t_0$ where $t_0$ is defined by $t_0=\arg\sinh (w_0^{-1/2})$.

Finally, by \eqref{F}
\begin{align*}
\phi_\lambda(e^t) \asymp \lambda^{-1} (
 -(\sinh t)^{-(\rho+\lambda)} + (\sinh t)^{\lambda-\rho})
 \asymp  \lambda^{-1}\, e^{(\lambda-\rho)t}
\end{align*}
what proves the conjecture in this case. 

\noindent{\bf Region IV:}  Suppose that  $t\ge c_1>0, \lambda \ge c_2>0$.

We use the formula \eqref{eq:shimeno}.   The uniform asymptotic approximation \cite[(3.76), p. 689]{Olde}, valid uniformly  for   large $|z|$ when $\lambda_0\to \infty$,  
implies the   following estimate when  $t\ge c_1>
0, \lambda_0 \ge c_2>0$,  for $c_1,c_2$ large enough: 
\begin{align*}
\;_2 F_ 1(a +\lambda_0 ,a- \lambda_0,
c
;-z) \asymp 
\lambda_0^{ 1/2-c} \zeta^ {c-2a-1/2} z^{ -a}
\left((z(1 + \zeta)^2)^{\lambda_0} +
\frac1{
(z(1 +\zeta)^2)^{\lambda_0}} \right),
\end{align*}
where $a=\rho/2$, $\lambda_0=\lambda/2$, $c=k_1+k_2+1/2$, $z=\sinh^2 t$
and $\zeta=(1+z^{-1})^{1/2}= \coth t$.
Consequently, 
$z(1 + \zeta)^2=(\sinh t + \cosh t)^2=e^{2t}$.
We get
\begin{align*}
\phi_\lambda(e^t)\asymp \lambda^{ -k} 
 \coth^ {k-\rho} t \;
 e^{ -\rho t}
(e^{\lambda t} +
e^{-\lambda t} )
\asymp \lambda^{ -k} 
 e^{ (\lambda-\rho) t}
\end{align*}
 which proves the proposition in this case.

The closure of the complement of the four regions we have discussed above in the set $\overline{\a^+}\times \overline{\a^+}$ is a compact set.  Given that $\phi_\lambda(e^X)$ and the proposed bounds are both continuous in $(\lambda,X)$ and nonzero (for $\phi_\lambda$, this can be easily seen from \eqref{Koor}), the result follows.
\end{proof}

\begin{remark}
Asymptotic expansions and approximations of the hypergeometric function $\;_2F_1$ for large values of parameters
is an important and active research topic research see \cite{Bateman}, \cite[Chapter 15.12]{NIST} and \cite{Olde} as a survey of results.

In \cite{Watson}, Watson gave an asymptotic expansion of the function
$\;_2F_1(a+\lambda, b-\lambda, c; (1-z)/2)$ for large $|\lambda|$,
resumed in  \cite[(17) p. 77]{Bateman}.
His result also yields the estimate  $\phi_\lambda(e^t)\asymp \lambda^{ -k} e^{ (\lambda-\rho) t}$
for large $\lambda$ and $z$ although the uniformity on his bound in $z$ is not clearly stated.
\end{remark}

\subsection{Proof of Conjecture \ref{conj:all} in a region for a symmetric space of noncompact type}

We first recall a result from a previous paper (\cite[Proposition 3.5]{PGPSTL}).

\begin{proposition}\label{XWY}
Let $\alpha_i$ be the simple roots and let $A_{\alpha_i}$ be such that $\langle X,A_{\alpha_i}\rangle=\alpha_i(X)$ for $X\in\a$.
Suppose 
$X\in\a^+$ and $w\in W\setminus\{id\}$. Then we have
\begin{align}\label{CL}
Y-w\,Y=\sum_{i=1}^r \,2\,\frac{a_i^w(Y)}{|\alpha_i|^2}\,A_{\alpha_i}
\end{align}
where $a_i^w$ is a linear combination of positive simple roots with non-negative integer coefficients for each $i$.
\end{proposition}

\begin{corollary}\label{C}
Suppose $X$, $\lambda\in\overline{\a^+}$.  Then there exists $M>0$ depending only on the Lie algebra structure such that for $H\in C(X)$, the convex hull of $W\cdot X$,
$\lambda(X)-M\,\max_{1\leq i,j\leq n}\,\{\alpha_i(\lambda )\,\alpha_j(X)\}\leq \lambda(H)\leq \lambda(X)$.
\end{corollary}
\begin{proof}
Since $\lambda$ is a linear function, it attains its maximum at the extremal points of $C(X)$ namely on $W\cdot X$.  The rest follows from the Proposition \ref{XWY}.
\end{proof}

The next result shows, using the well known estimates for $\phi_0(e^X)$,  that Conjecture \ref{conj:all} holds in a region of the variables $(\lambda,X)$ for a symmetric space of noncompact type.

\begin{proposition}
Let  $X$, $\lambda\in\overline{\a^+}$ and suppose $\alpha(\lambda)\,\alpha(X)\leq C$ for all $\alpha \in \Sigma^+$. Then there exists $M>0$ depending only on the Lie algebra structure 
\begin{align*}
e^{-M\,C}\,e^{\lambda(X)}\,\phi_0(e^X)\leq \phi_\lambda(e^X)\leq e^{\lambda(X)}\,\phi_0(e^X)
\end{align*}
\end{proposition}  

\begin{proof}
We consider every classical or exceptional Lie algebras.  According to \cite[Plates I--IX]{Bourbaki}, every root system has a highest root of the form 
$\gamma=\sum_{i=1}^r\,n_i\,\alpha_i$ where $n_i\geq 1$ for each $i$.  The condition $\gamma(\lambda)\,\gamma(X)\leq C$
then implies that $\max_{1\leq i,j\leq n}\,\{\alpha_i(\lambda )\,\alpha_j(X)\}\leq  C$.

Now, 
\begin{align*}
\phi_\lambda(e^X)=\int_K\,e^{(\lambda-\rho)(H(e^X\,k))}\,dk=\int_K\,e^{\lambda(H(e^X\,k))}\,e^{-\rho(H(e^X\,k))}\,dk.
\end{align*}

Noting that $\{H(e^X\,k)\colon k\in K\}=C(X)$, the result follows from Corollary \ref{C}.
\end{proof}

\section{Acknowledgements}
We thank  A. B. Olde Daalhuis and A. Nowak for advice on asymptotics of the hypergeometric function.

We are grateful to the grants IEA CNRS: Analyse li\' ee aux racines et applications  2021--2022 and MIR Universit\'e d'Angers ``Sym\'etries'' for their support of this research.

\end{document}